\documentclass[12pt]{article}
\usepackage{graphicx,amssymb,amsmath,amsthm}
\usepackage{hyperref,url}
\usepackage{verbatim}
\title{Filling a triangulation of the 2-sphere}
\author{Peter Doyle \and Matthew Ellison \and Zili Wang}
\date{Version 3.0 dated 17 May 2025}
\newtheorem*{hypo}{Hypothesis}

\newtheorem*{conjecture}{Conjecture}

\newcommand{\beginfigurepos}{\begin{figure}[htbp]}
\newcommand{\fig}[3]{
\beginfigurepos
\includegraphics[width=370pt]{figures/#1.pdf}
\caption{#3}
\label{#2}
\end{figure}
}
\newcommand{\figsize}[4]{
\beginfigurepos
\centerline{
\includegraphics[width=#1]{figures/#2.pdf}
}
\caption{#4}
\label{#3}
\end{figure}
}

\def\del{\partial}

\newcommand{\ceiling}[1]{\left \lceil #1 \right \rceil}

\makeatletter
\newcommand{\xRrightarrow}[2][]{\ext@arrow 0359\Rrightarrowfill@{#1}{#2}}
\newcommand{\Rrightarrowfill@}{\arrowfill@\equiv\equiv\Rrightarrow}
\newcommand{\xLleftarrow}[2][]{\ext@arrow 3095\Lleftarrowfill@{#1}{#2}}
\newcommand{\Lleftarrowfill@}{\arrowfill@\Lleftarrow\equiv\equiv}
\newcommand{\xLleftRrightarrow}[2][]{\ext@arrow 3399\LleftRrightarrowfill@{#1}{#2}}
\newcommand{\LleftRrightarrowfill@}{\arrowfill@\Lleftarrow\equiv\Rrightarrow}
\makeatother

\newcommand{\kw}{\mathbf}

\newcommand{\tetvol}{\kw{tetvol}}
\newcommand{\Qvol}{\kw{Qvol}}

\newcommand{\vol}{\kw{vol}}
\newcommand{\flipdist}{\kw{flipdist}}

\newcommand{\maxdeg}{\kw{maxdeg}}

\newcommand{\vsa}{\kw{vsa}}

\begin{document}

\maketitle

\begin{abstract}
Define the tet-volume of a triangulation of the 2-sphere to be
the minimum number of tetrahedra in a 3-complex
of which it is the boundary, and
let $d(v)$ be the maximum tet-volume for $v$-vertex triangulations.
In 1986 Sleator, Tarjan, and Thurston (STT)
proved that $d(v) = 2v-10$ holds for large $v$,
and conjectured that it holds for all $v \geq 13$.
Their proof used ideal hyperbolic polyhedra of large volume.
They suggested using more general notions of volume instead.
In work that was all but lost, Mathieu and Thurston
used this approach to outline a combinatorial proof of the STT
asymptotic result.
Here we use a much simplified version of their approach
to prove the full conjecture.
\begin{comment}
This implies STT's weaker conjecture,
proven by Pournin in 2014,
characterizing the maximum rotation distance between trees.
\end{comment}
\end{abstract}
\centerline{\emph{For Bill}}

\section{Cover letter}

In their famous paper
`Rotation distance, triangulations, and hyperbolic geometry'
\cite{stt:jams},
Sleator, Tarjan, and Thurston (STT)
pointed out the equivalence between the rotation distance between
two $v-2$-vertex binary trees and the flip distance
between two triangulations $\sigma_1,\sigma_2$ of a $v$-gon,
and got lower bounds for $\flipdist(\sigma_1,\sigma_2)$
by using hyperbolic geometry to derive lower bounds for the number
of tetrahedra needed to complete the triangulation
$\sigma=\sigma_1-\sigma_2$ of the sphere $S^2$ to a triangulation of 
the ball $B^3$.

This paper has been cited over 500 times.
Most citers were likely looking only for what the editors of
Thurston's collected works
\cite[p. 331]{thurston:works3}
call the `astonishingly elementary observation'
that rotations of a binary tree correspond to edge flips of a triangulated
polygon, or the `equally elementary argument' that for $v \geq 13$,
$\flipdist(\sigma_1,\sigma_2) \leq 2v-10$.
But some, at least, came to marvel at the `deep and original
argument involving volumes of hyperbolic 3-manifolds'
that `solved using transcendental methods a purely combinatorially stated
problem for which good asymptotic estimates were not even known.'

STT didn't fully solve the flip distance problem
that had inspired them to look at triangulations of the sphere,
because their proof worked only for $v$ sufficiently large.
Eventually Pournin
\cite{pournin:diameter}
solved the flip distance problem,
but that didn't settle STT's conjecture that
for all $v \geq 13$
there is a $v$-vertex simplicial triangulation $\sigma$ of $S^2$
that cannot be extended to a simplicial triangulation of $B^3$
with fewer than $2v-10$ $3$-simplices.

In an addendum to their paper, STT
\cite[p. 697]{stt:jams}
suggested
using a
purely combinatorial alternative
to considerations of hyperbolic volume.
Instead of the volume of a hyperbolic
realization of $\sigma$, look at $\Qvol(\sigma)$,
the minimum $L_1$-norm of a
$3$-chain in the simplex $\Delta^{v-1}$ with boundary $\sigma$.
This is a rational number, the value of a linear program.
Mathieu and Thurston (MT) explored this approach in an extended abstract
that outlined an alternative proof of the asymptotic lower
bound $2v-10$ obtained by STT.
This abstract was submitted to and rejected by STOC in 1992;
the authors turned to other things;
and the work was all but lost.

In this paper we take up this abandoned approach,
and use it to prove STT's conjecture that for all $v \geq 13$
there is a $v$-vertex simplicial triangulation of $S^2$
that cannot be extended to a simplicial triangulation of $B^3$
with fewer than $2v-10$ $3$-simplices.

\section{Summary}

A \emph{triation of the sphere} is
an oriented simplicial 2-complex $\sigma$
whose carrier is the $2$-sphere.
Regard $\sigma$ as 
a subcomplex of the $(v-1)$-simplex $\Delta^{v-1}$,
where $v$ is the number of vertices.
A \emph{tetration} of $\sigma$ is 
an oriented 3-subcomplex $\tau$ of $\Delta^{v-1}$ with boundary
$\del \tau = \sigma$.
Define the \emph{tet-volume} $\tetvol(\sigma)$ to be the minimum number of 
`tets' (tetrahedra)
in a tetration.
Let $d(v)$ be the maximum of $\tetvol(\sigma)$ over all $v$-vertex triations
of the $2$-sphere.

In 1986 Sleator, Tarjan, and Thurston (STT) observed
\cite{stt:stoc,stt:jams}
that
if $v \geq 13$ then $d(v) \leq 2v-10$,
as follows.
From Euler, a triation with $v$ vertices has
$2v-4$ faces.
For any choice of vertex $a$ we get a tetration by coning from $a$.
(Take the union of all tets $[abcd]$ where $[bcd]$ is an oriented
face of $\sigma$ disjoint from $a$.)
This shows
\[
\tetvol(\sigma) \leq 2v-4-\deg(a)
,
\]
so choosing $a$ to achieve the maximum degree $\maxdeg(\sigma)$
we have
\[
\tetvol(\sigma) \leq 2v-4-\maxdeg(\sigma)
.
\]
Assuming $v \geq 13$ we have $\maxdeg(\sigma) \geq 6$, yielding
\[
\tetvol(\sigma) \leq 2 v-10
.
\]

STT conjectured:
\begin{conjecture}[tetvol]
\[
d(v)=2v-10
,\; v\geq 13
.
\]
\end{conjecture}
\noindent
They proved this with $13$ replaced by some unspecified constant,
by considering triations arising as the boundary of certain ideal hyperbolic
polyhedra with large volume.
Because there is an upper bound for the volume of a hyperbolic tet,
to tetrate a hyperbolic polyhedron with large volume requires a large number
of tets.
The lower bound they get is only $2v-O(\log(v))$,
but having got themselves into the ballpark,
they proceed to pull off an inside-the-park home run,
showing that not only is $\tetvol=2v-10$,
but any optimal tetration comes from coning.

In an addendum to their paper,
STT \cite[p. 697]{stt:jams}
suggested that instead of using hyperbolic volume,
lower bounds could be proven purely combinatorially.
Mathieu and Thurston (MT) pursued this approach in
\cite{mt}.
(See section \ref{mt}.)
We take this same approach here,
and prove the tetvol conjecture
by examining a very simple class of triations, obtained
by truncating skinny cylindrical quotients of the Eisenstein lattice.

The work of STT on tetrations
was motivated by the problem of
finding the maximum rotation distance between trees,
or equivalently,
the maximum flip distance $d'(v)$ between
two triangulations of a $v$-gon.
(See \cite{stt:jams} for definitions and discussion.)
The tetvol conjecture implies the associated conjecture that
$d'(v)=2v-10$ for $v \geq 13$.
(See section \ref{flipdist}.)
In 2014 Pournin
\cite{pournin:diameter}
proved this associated conjecture directly,
without reference to tet-volume,
by producing pairs
that he could show have 
flip distance $2v-10$.
Pournin's theorem doesn't settle the tetvol conjecture,
because there may be a gap between tet-volume and flip distance---
see section \ref{flipdist} below.
The examples described here provide myriad pairs
maximizing flip distance, with Pournin's pairs among them.

Here's an outline of the paper.
First we'll introduce the key idea of a volume potential,
as applied to the icosahedron, a toy case that shows the basic idea.
Then we'll apply the method to our skinny triations
to settle the tetvol conjecture.
We'll then discuss our debt to MT;
the connection to linear programming;
related constructions of triations;
implications for flip distance;
and where we go from here.

\section{Warming up with the icosahedron}
\newcommand{\icos}{\mathbf{icos}}

We'll begin with the icosahedron
($v=12$; $f=20$; $\maxdeg=5$).
Coning from a vertex
yields a tetration with $20-5= 15=2v-9$ tets,
so $\tetvol(\icos) \leq 15$.
Now let's show that $\tetvol(\icos) \geq 15$.

We will call a function $\rho(ABC)$ on ordered triples
of vertices $ABC$ satisfying
\[
\rho(ABC)=\rho(BCA) = - \rho(ACB)
\]
a \emph{volume potential}.
For any permutation $XYZ$ we get from this that
$\rho(XYZ)= \pm\rho(ABC)$, depending on the sign of the permutation.
(In standard language, $\rho$ is a 2-cochain.)

Associated to a volume potential is its
\emph{volume form}
$\vol_\rho(ABCD)$, the function on ordered triples $ABCD$ given by
\[
\vol_\rho(ABCD) = \rho(BCD)-\rho(ACD)+\rho(ABD)-\rho(ABC)
.
\]
For any permutation $XYZW$ of $ABCD$ we have $\vol_\rho(XYZW) = \pm\vol_\rho(ABCD)$.
(In standard language,
$\vol_\rho$ is 3-cocycle, the coboundary of $\rho$.)

Let $\rho(\icos)$ be the sum of $\rho(ABC)$ over the faces $ABC$ of $\icos$.
(Here and hereafter, by `faces' we mean properly oriented faces.)
For $\tau$ a tetration of $\icos$
let $\vol_\rho(\tau)$ 
be the sum of $\vol_\rho(ABCD)$ over the tets of $\tau$.
The key fact we need is that
\[
\vol_\rho(\tau) = \rho(\icos)
.
\]
This is Stokes's theorem in combinatorial form;
it's true because matching faces of the tets of $\tau$ make
contributions of opposite sign to $\vol_\rho(\tau)$,
so that after cancellation only the contributions from
faces of $\icos$ remain.

Call the volume potential $\rho$ \emph{good} if the volume form
$\vol_\rho$ assigns all tets $ABCD$ volume at most $1$:
\[
\vol_\rho(ABCD) \leq 1
.
\]
Observe that this implies that
\[
|\vol_\rho(ABCD)| \leq 1,
\]
since
\[
\vol_\rho(ABCD)= -\vol_\rho(ABDC) \geq -1
.
\]

For a good volume potential $\rho$ the number $|\tau|$ of tets of $\tau$
satisfies
\[
|\tau| \geq \vol_\rho(\tau)  = \rho(\icos)
.
\]
So to prove that $\tetvol(\icos)=15$, we just need to find
a good volume potential $\rho$ with $\rho(\icos)=15$.
(Or at least with $\rho(\icos)>14$, because we can always round a non-integral
lower bound up:
see section \ref{LP}.)

Let's look for a volume potential $\rho$ that is
invariant under orientation-preserving
symmetries of $\icos$.
(This won't hold us back: If there's any $\rho$ at all, we can
average to get a symmetrical one.)

To define $\rho$,
we need to prescribe a value for every triple $ABC$ of
distinct vertices.
We distinguish three cases.
\begin{itemize}
\item
$ABC$ or its orientation-reversal $ACB$ is a face of $\icos$.
Since we want $\rho$ to be symmetric and $\rho(\icos)=15$,
we must take $\rho(ABC)=3/4$ if $ABC$ is a face,
which makes $\rho(ABC)=-3/4$ if $ACB$ is a face.
\item
$\{B,C\}$ is an edge of $\icos$, but not both $\{A,B\}$ and $\{A,C\}$.
This is the crucial case.
By symmetry, we can assume that $A$ is any fixed vertex of $\icos$.
Figure \ref{icos}
shows the values for a fixed choice of $A$ by means of a flow $\phi$
along the edges of the dual graph.
Along the dual edge $BC^\perp$ clockwise from $BC$ the flow rate
is $\phi(BC^\perp) = \rho(ABC)$.
\figsize{.5\linewidth}{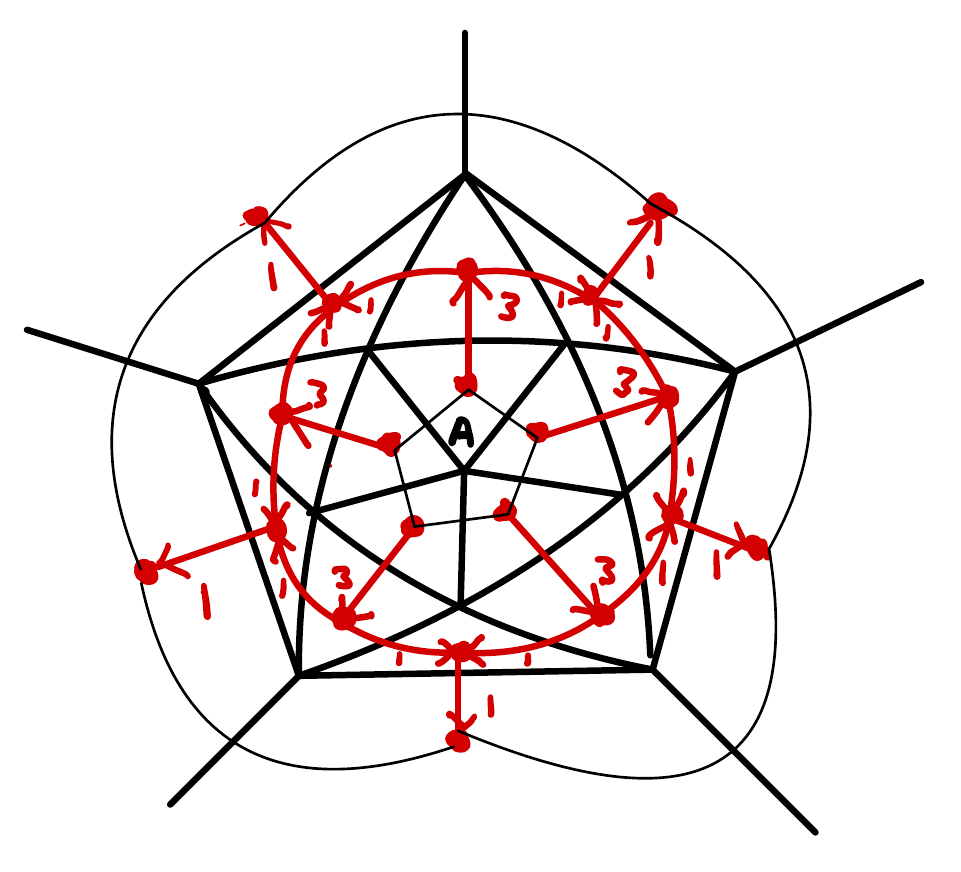}{icos}{
Volume potential flow for $\icos$.
Flow rates are one quarter of what is shown.
To an edge $BC$ we assign a flow of rate $\rho(ABC)$
along the dual edge $BC^\perp$.
This flow is not conservative:
There is net flow $3/4$ out of each of the five faces incident with $A$,
and net flow $1/4$ into each of the fifteen remaining faces.
}
\item
$ABC$ involves no edge of $\icos$.
Here we'll take $\rho(ABC)=0$.
Up to symmetry of $\icos$
there is just one possibility for the unoriented triangle with vertices
$\{A,B,C\}$---or
two possibilities if you don't allow orientation-reversing symmetries---but
for this proof we don't need to check this.
\end{itemize}

Now we want to check that $\rho$ is good, i.e.,
that $\vol_\rho(ABCD) \leq 1$ for any 4-tuple $ABCD$.
Suppose $ABCD$ contains a face of $\icos$.
We may assume that $BCD$ is a properly oriented face,
with $A$ having been moved to the standard position.
We have cooked up the flow
$\phi$
so that the net flow into any face of $\icos$
not incident with $A$ is $1/4$.
In terms of $\rho$ the net flow into $BCD$ is
\begin{eqnarray*}
1/4
&=&
-\phi(BC^\perp)-\phi(CD^\perp)-\phi(DB^\perp)
\\&=&
-\rho(ABC)-\rho(ACD)-\rho(ADB)
\\&=&
-\rho(ABC)-\rho(ACD)+\rho(ABD)
.
\end{eqnarray*}
Since $\rho(BCD) = 3/4$,
this gives us
\[
\vol_\rho(ABCD) = \rho(BCD)-\rho(ACD)+\rho(ABD)-\rho(ABC) = 3/4+1/4=1
.
\]

We're now close to having shown that $\rho$ is good.
We just need to check the case where the tet
$ABCD$ contains no face of $\icos$,
whether properly or improperly oriented.
But for such tets the four terms of $\vol_\rho(ABCD)$
all have absolute value at most $1/4$,
so $|\vol_\rho(ABCD)| \leq 1$.

Having verified that $\rho$ is good, with $\rho(\icos)=15$, we're done.

\section{Proof of the tetvol conjecture} \label{done}
To prove the tetvol conjecture,
all we need is a sequence of triations
$T_v$ and good volume potentials $\rho_v$
with $\rho_v(T_v) = 2v-10$ for all $v\geq 13$.
These are defined and checked by the code in Figures
\ref{Tcode} and \ref{rhocode}, which together constitute
a proof of the tetvol conjecture.
Or rather, a `verification'.
Where's the proof?
\begin{figure}[tbp]
\begin{verbatim}
#!/usr/bin/python3

from rho import rho # import the volume potential rho

# define phyllohedra T_v

def Tbase(v):
  return [(0,5,4),(0,4,3),(0,3,2),(0,2,1)]
def Tsides(v):
  pointup=list((k,k+1,k+6) for k in range(v-6))
  pointdown=list((k,k+6,k+5) for k in range(v-6))
  return pointup+pointdown
def Tlid(v):
  return [(v-1,v-6,v-5),(v-1,v-5,v-4),(v-1,v-4,v-3),(v-1,v-3,v-2)]
def T(v): return Tbase(v)+Tsides(v)+Tlid(v)

# we need only check rho(T_v)=2v-10 for two values of v

def rhovol(poly): return sum(rho(abc) for abc in poly)
assert rhovol(T(13))==2*13-10 and rhovol(T(14))==2*14-10
\end{verbatim}
\caption{Code to define $T_v$ and check that $\rho(T_v)=2v-10$.}
\label{Tcode}
\end{figure}

To understand the family $T_v$,
let's start by looking at some pictures:
Figures \ref{t23} and \ref{t51}.
(Better yet, build some physical models!)
The vertices of $T_v$ are labeled $0,\ldots,{v-1}$ spiraling up from the bottom.
Cutting along edges $01$ and $(v-2)(v-1)$ we get a topological
cylinder which unwraps to give a diagram like those shown
in Figure \ref{unwrap} for $v=13$ and $v=14$.
\figsize{.8\linewidth}{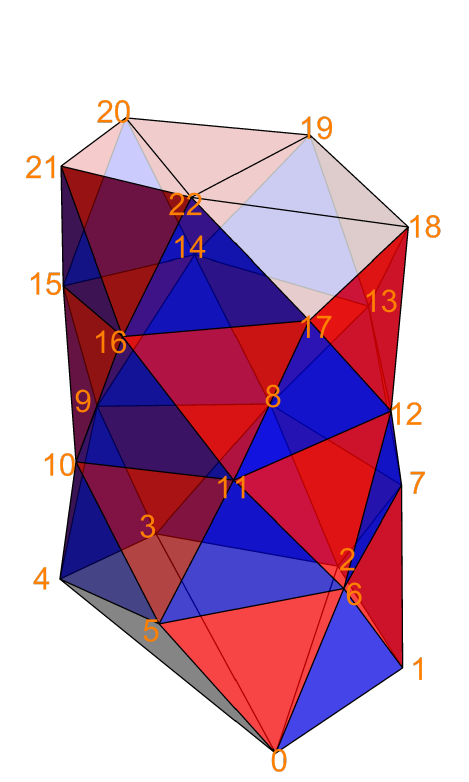}{t23}{Triation $T_{23}$}
\figsize{\linewidth}{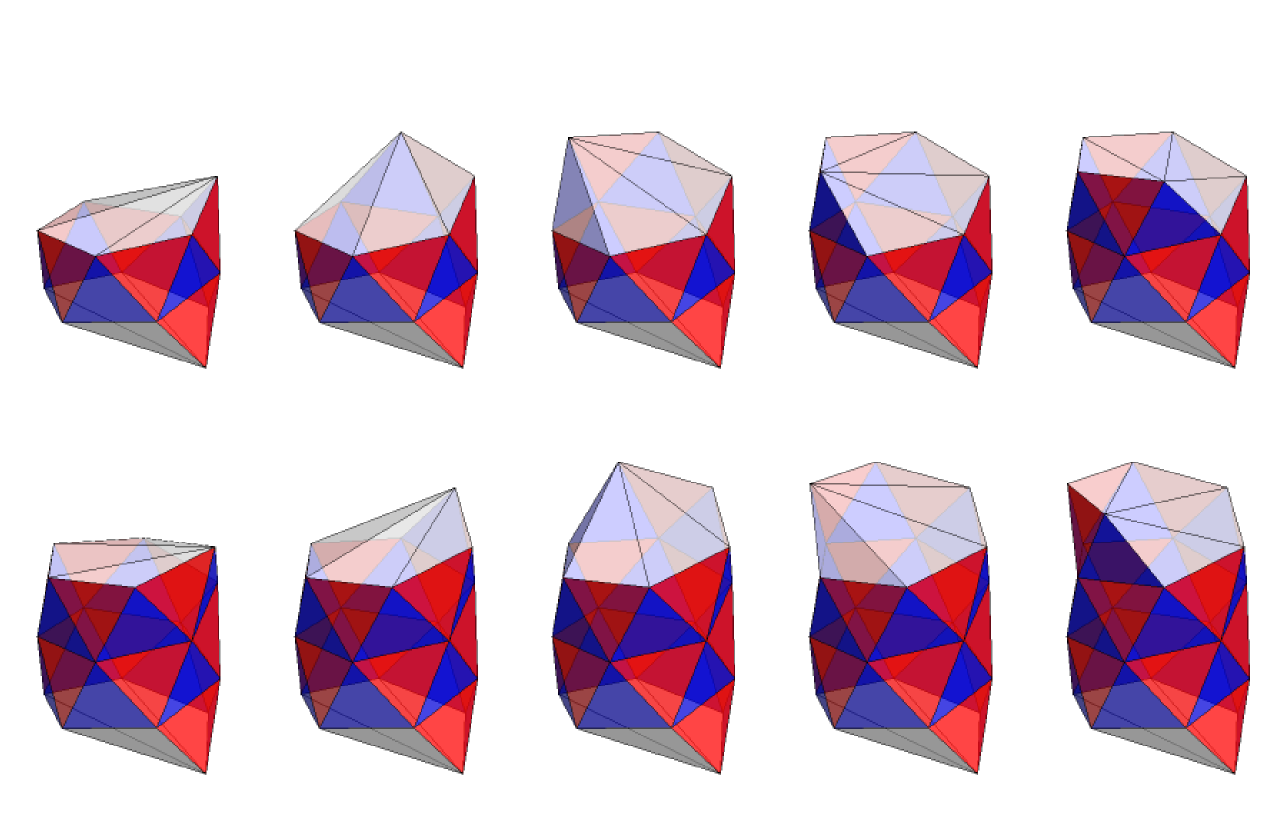}{t51}{Triations $T_{13}$ through $T_{22}$}
\figsize{.8\linewidth}{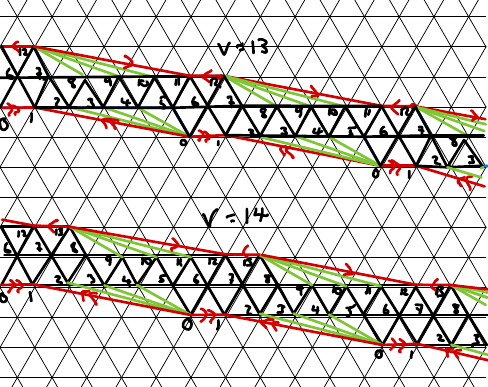}{unwrap}{$T_{13}$ and $T_{14}$ unwrapped}

These triations $T_v$ are examples of what we call `phyllohedra'.
They are obtained as follows.
\newcommand{\Eis}{\mathbf{Eis}}
Associated to the Eisenstein integers
$\Eis = \{u + v \omega$\}, $\omega = \exp(i \tau/3)$
is a triation of the plane with six triangles meeting per vertex.
This triation descends to the quotient cylinder
\[
\Phi_{a,b} = \Eis / ((a-b\omega)\mathbf{Z})
.
\]
Wrapping $\Phi_{a,b}$ around the $z$-axis,
we get pictures
like those that arise in phyllotaxis.
(See Figure \ref{phyllocyl}.)
This inspires us to call $\Phi_{a,b}$ \emph{the $(a,b)$-phyllocylinder}.

\figsize{\linewidth}{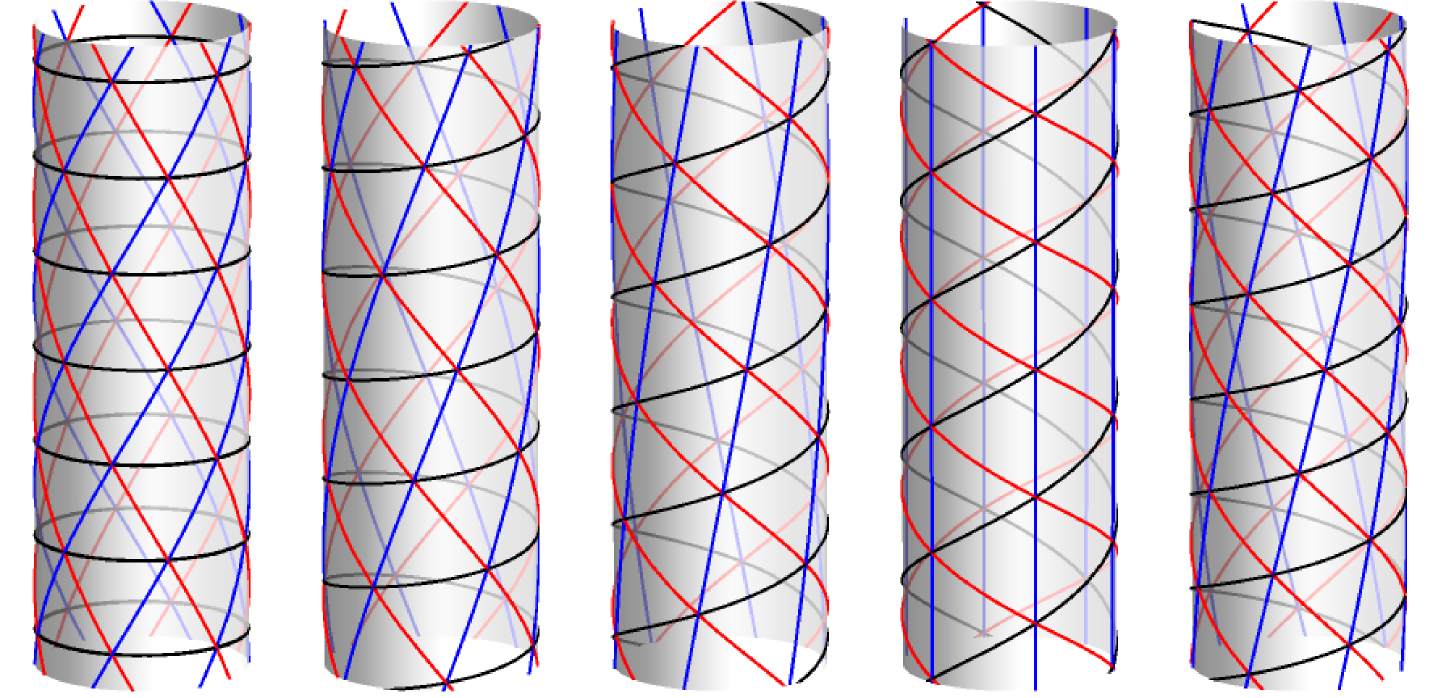}{phyllocyl}{
The phyllocylinders $\Phi_{6,0},\Phi_{5,1},\Phi_{4,2},\Phi_{4,2},\Phi_{5,2}$.
Unless $a$ or $b$ vanishes,
$\Phi_{a,b}$ has
$b$ spirals in direction $1$ (black);
$a$ spirals in direction $\omega$ (red);
$a+b$ spirals in direction $1+\omega$ (blue).
}

We can truncate a phyllocylinder by taking a subset of the vertices
and the triangles they inherit
from $\Phi_{a,b}$,
together with some extra edges and faces
to cap off the bottom and top.
We call these finite triations `phyllohedra', a loose term
whose precise meaning will depend on what kinds of truncation
and capping you allow.

Our triations $T_v$ are $(5,1)$-phylohedra, obtained
by truncating $\Phi_{5,1}$ and then capping in the most natural way.

\beginfigurepos
\begin{verbatim}
#!/usr/bin/python3

def rho(abc):
  (a,b,c)=abc

  # make sure a<=b<=c
  if a>b: return -rho((b,a,c))
  if b>c: return -rho((a,c,b))

  # rho depends only on the gaps
  (x,y)=(b-a,c-b)
  if x>5:
    return 0
  elif 3<=x<=5 and x+y>=6:
    return -1
  elif x==2 and y==1:
    return 1
  elif x==1 and b>=1:
    return 1
  else:
    return 0

def vol(abcd):
   (a,b,c,d)=abcd
   return rho((b,c,d))-rho((a,c,d))+rho((a,b,d))-rho((a,b,c))

# check goodness of rho for 4-tuples of integers between 0 and n-1
# n=18 should do it, but let's go way overboard

n=36
tuples=list((a,b,c,d) for a in range(n) for b in range(n)
                      for c in range(n) for d in range(n))
vols=list(vol(abcd) for abcd in tuples)

# for this rho the only values we should see are 0,1,-1
assert set(vols)=={0,1,-1}

\end{verbatim}
\caption{Code to define $\rho$ and check that it is good.}
\label{rhocode}
\end{figure}
To accompany our $T_v$, we need volume potentials
$\rho_v$.
We will take these to be restrictions of a single translation-invariant
volume potential $\rho$ defined on the infinite cylinder $\Phi_{5,1}$.
This is possible because the vertices of $T_v$ are a subset of the vertices
of $\Phi_{5,1}$.
On any $\Phi_{a,1}$ the vertices are nicely indexed by integers,
so we can think of $\rho$ as defined for triples of integers.
Translation invariance means that
$\rho((a,b,c))=\rho((a+k,b+k,c+k))$.

Figure \ref{rhocode} shows code to compute $\rho$,
and check the volume condition for all 4-tuples of
integers between $0$ and $n-1$.
Because of the way $\rho((a,b,c))$ depends only on the gaps
between $a,b,c$,
and treats gaps that are $6$ or bigger as equal,
taking $n=18$ should cover all possible cases;
we take $n=36$ in case $n=18$ is off by one
(or two, or three, or four,\ldots).

To check the goodness of $\rho$ by hand,
we can look at
the volume potential flow,
as we did for the icosahedron.
By the translation invariance of $\rho$,
here again we need only a single picture:
Figure \ref{x51}.
This time the inflow vanishes for triangles other than those containing
the reference vertex $A$, here represented by a black dot.
This makes $\vol_\rho(ABCD)=1$ whenever $BCD$ is a face of $\Phi_{5,1}$.
\figsize{.8\linewidth}{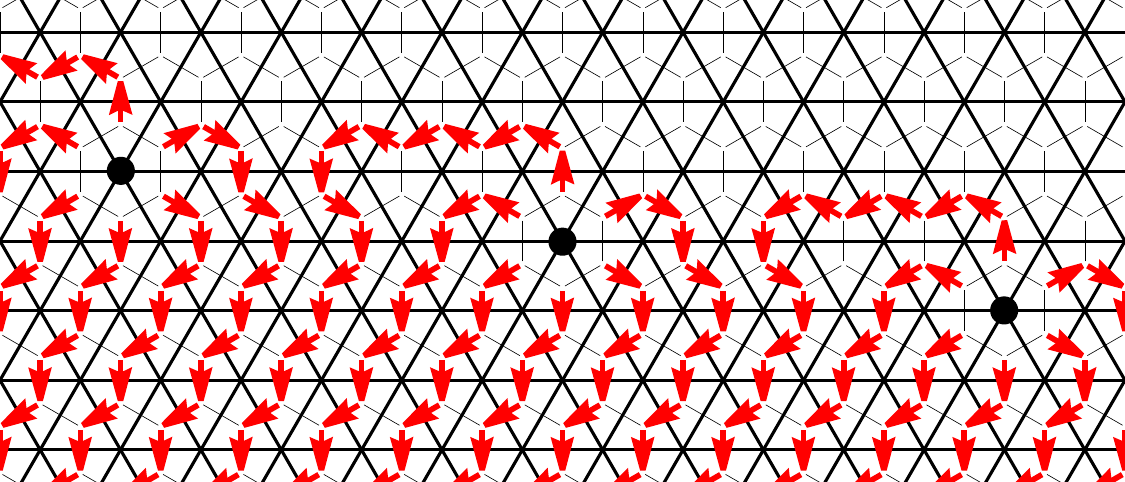}{x51}{Volume potential flow for
$\Phi_{5,1}$. Reference vertex $A$ is the black dot.
Flow along red arrows has rate 1. The net flow into triangles not incident
with $A$ is $0$.}
This takes us a long way toward showing that $\rho$ is good.
We still have to deal with tets $ABCD$ not involving any face of $\icos$,
and hence with values $\rho$ of triangles not containing any edge.
We don't know any really clever way to check these `big' tets.
Fortunately lots of triangles have $\rho = 0$,
which makes things easier.
We can either roll up our sleeves and get to work,
or decide to trust the computer on this.

To complete the proof,
we must check that $\rho(T_v) = 2v-10$.
Of the $2v-4$ `side' faces of $T_v$, all but eight are faces of $\Phi_{5,1}$,
and thus get weight $1$, which gets us up to $2v-12$.
We just need to check that the net contribution of the base and lid
together give us the extra $2$ we need.
Since these faces all contain an edge of $\Phi_{5,1}$
(either $01$ or $(v-2)(v-1)$)
the information needed to check this is there in the flow diagram
\ref{x51}.
As an alternative, or as a check on our work, we can observe that just have to
check a single value of $v$ to nail down the constant term, and the
code in Figure \ref{Tcode} has done this for us.

This completes the proof.

We should emphasize that this $\rho$ is not canonical,
not even at the level of the associated volume form $\vol_\rho$.
\begin{comment}
This $\rho$
gives a tight lower bound for $\tetvol$
for this particular family of $(5,1)$-phyllohedra.
For other families you may need to choose a different $\rho$---perhaps
one that is not translation-invariant.
There might not even be a single $\rho$ that works for all the triations
in your family.
For all we know, there might not be any $\rho$ at all:
See section \ref{LP}.
\end{comment}
It's particularly nice in that it takes
only values $0,1,-1$.
In fact you can get the exact value of $\tetvol$ for small triations
with such `binary' volume potentials.
But once you get up to $v=19$ or so, you find triations with non-integral
$\tetvol$, and then the jig is up.

\section{Volume potentials} \label{mt}

The volume potential method we're using here was proposed
(in an equivalent form) by STT;
In the early 1990s Mathieu and Thurston (MT)
explored this approach, and introduced the flows
we've used to encode values of the volume potential.
As in STT, MT did not prove the full tetvol conjecture,
just the version for $v$ large enough.
Their paper was not published, and eventually was all but lost.
The work described here was based on a fragment of the paper
gleaned from a garbled fax,
and the email from Bill Thurston reproduced in Figure \ref{billmail}.
%and some email from Bill Thurston.
Recently a complete copy of the paper was found
and posted to the arxiv
\cite{mt}.

\begin{figure}[p]
{\fontsize{9}{8}
\begin{verbatim}
From doyle Thu Nov 14 23:14:15 2002
To: wpt
Subject: rotations

Bill,

I understand that there was a draft or a preprint related to the attached
abstract.  If you can send anything (e.g. a tex source) I'd love to see it.

Peter
-----------

Claire Kenyon (ENS-Lyon and William Thurston, MSRI)

Rotation distance between binary trees: hyperbolic geometry vs. max-flow
min-cut.

The maximum number of rotations needed to go from one given binary tree
with $n$ nodes to another is exactly $2n-6$ when $n$ is large enough.
We first sketch Sleator, Tarjan, and Thurston's original proof of this
theorem, which involves hyperbolic geometry volume arguments, the
Riemann mapping theorem, approximate calculations of integrals and
an induction argument. We then present an alternate, elementary proof,
based on the max-flow min-cut theorem. Finally, we compare the two
proofs and show how they are essentially two versions of the same
proof, by relating successively hyperbolic volume to cocycles to
linear programming to amortized analysis to flow problems.

-----

From wpthurston@mac.com  Fri Nov 15 01:25:13 2002
Date: Thu, 14 Nov 2002 22:23:40 -0800
Subject: Re: rotations
From: wpthurston@mac.com
To: "Peter G. Doyle" <doyle@hilbert.dartmouth.edu>

Hi Peter,
Yes, there was a draft, but in all my moving around I don't have a copy any more.
I should try to get my own copy from Claire.  We kind of dropped it when it 
was rejected from STOC.

The idea was, given a triangulation that appears to be maximal, attempt to
construct an L^infty 3-cocycle on \Delta^{v-1}  (if there are v vertices) that
takes maximal value on all the tetrahedra in your triangulation.  Of course, hyperbolic
volume (given an immersion of the polyhedron into H^3) gives a cocycle that
works well enough in many cases --- but it's not quite the best.
I don't remember all the details, although I could reconstruct them.  I think
it turned out that the important values to work out were when 2 or more vertices
of the tetrahedron are connected by an edge; these could be done using
many instances of a max-flow min-cut process: I think, one instance for each
possible location for the pair of non-adjacent vertices (it was like a flow
from one of these vertices to the other). When three vertices are all mutually
connected by edges (i.e. the tetrahedron has a face on the bounding sphere)
I think there is some formula that we just wrote down.  I think the cocyle could
assign 0 to many of the other tetrahedra, the ones with 4 non-adjacent vertices.

Of course a cocycle like this is really the dual to the L^1 3-chain having boundary
the given triangulation. There might or might not be a geometric triangulation
realizing this minimum, but it seemed to work out in lots of cases, including
explicit examples for every value of the number of vertices such as in the
original paper.
        Bill
\end{verbatim}
}
\caption{The word from Bill}
\label{billmail}
\end{figure}

\begin{comment}
The work here did not follow directly in the footsteps of MT---but
it might as well have.
Good ideas that suggested themselves after exploring a raft of
not-as-good ideas show up in MT:
`So that's why they did that.'
\end{comment}

The MT paper was an `extended abstract'; It omitted certain details
meant to be covered in the `full paper' to follow.
We haven't tried to fill in the details, but
based on our own experience with this method,
we have no reason to doubt
that their work was essentially correct.

The difference between our approach and MT is that we deal with
simple examples where we can produce an explicit volume potential.
MT dealt with more complex examples, and used max-flow min-cut to find
the volume potentials.
We expect that their method will generalize in a way that ours will not.

\section{Linear programming} \label{LP}
Pick the volume potential $\rho$ so as to maximize $\vol_\rho(\sigma)$,
and call the maximum value $\Qvol(\sigma)$.
STT emphasized that this is a linear programming problem:
it's dual to the problem of minimally
tetrating $\sigma$ with fractional tets
allowed.

The value of $\Qvol$ is not always an integer.
(See section \ref{more}.)
When it isn't, we can round up:
\[
\tetvol(\sigma)
\geq
\ceiling{\Qvol(\sigma)}
.
\]
This suggests the null hypothesis:
\begin{hypo}[$\Qvol$]
\[
\tetvol(\sigma) =
\ceiling{\Qvol(\sigma)}
.
\]
\end{hypo}
This is false in general:
Ellison
\cite{ellison:hypo}
gives examples showing that the gap
$\tetvol(\sigma)-\ceiling{\Qvol(\sigma)}$
can be arbitrarily large.
These examples have $\maxdeg > 6$.
As far as we know, the hypothesis may still hold when $\maxdeg \leq 6$,
in which case
$\Qvol$ would still be a reliable way to detect triations
achieving the upper bound for $\tetvol$.
\begin{comment}
The problem is that things
only start to get interesting when you get to $v=20$ or so,
and we have difficulty computing $\tetvol$ or even $\Qvol$
for $v$ much bigger than this.
\end{comment}

\section{More about phyllohedra} \label{more}
Let's quickly indicate what happens when we truncate other phyllocylinders,
keeping the details for another day.

The phyllocylinder $\Phi_{a,b}$ has combinatorial girth $a+b$.
Along with $\Phi_{5,1}$, the other phyllocylinders of girth 6 are
$\Phi_{6,0}$, $\Phi_{4,2}$, and $\Phi_{3,3}$.
$\Phi_{6,0}$ is an infinite stack of hexagonal antiprisms.
Like $\Phi_{5,1}$, $\Phi_{6,0}$ has volume-to-surface-area ratio $1$,
meaning that it has a volume potential $\rho$ taking value $1$ to each of
its triangles:
$\vsa(\Phi_{6,0})=1$.
Figure \ref{x60} shows the associated flow.
(As usual it falls short in that it
doesn't indicate values for triangles not containing an edge.)
\begin{comment}
Truncating and capping $\Phi_{6,0}$
can yield examples proving the tetvol conjecture,
only now we need six families, as members of a family of $(a,b)$-phyllohedra
have $v$ congruent mod $\gcd(a,b)$.
\end{comment}
\figsize{.8\linewidth}{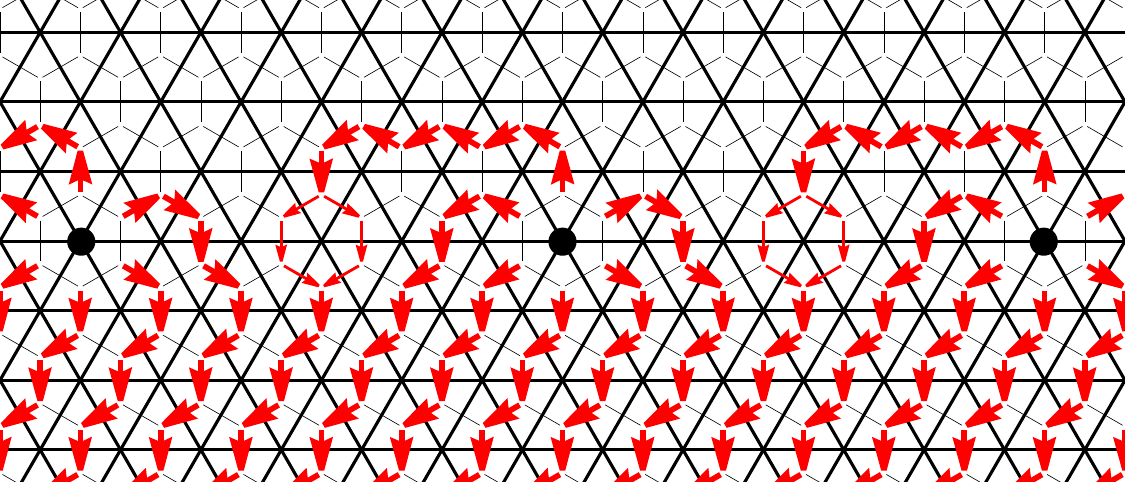}{x60}{Volume potential flow for $\Phi_{6,0}$.
Flow rate is $1$ along thick arrows, $1/2$ along thin arrows.}

By contrast,
\[
\vsa(\Phi_{4,2})=\frac{31}{32}
.
\]
Long $(4,2)$-phyllohedra have $\tetvol$ asymptotically equal to 
$\frac{31}{16} v$.
Truncating and capping in the most natural way
(see Figure \ref{icosplus}),
we get a family $U_v$ ($v$ even)
for which $\Qvol$ is not always an integer
--- but we still have $\ceiling{\Qvol}=\tetvol$.
\figsize{\linewidth}{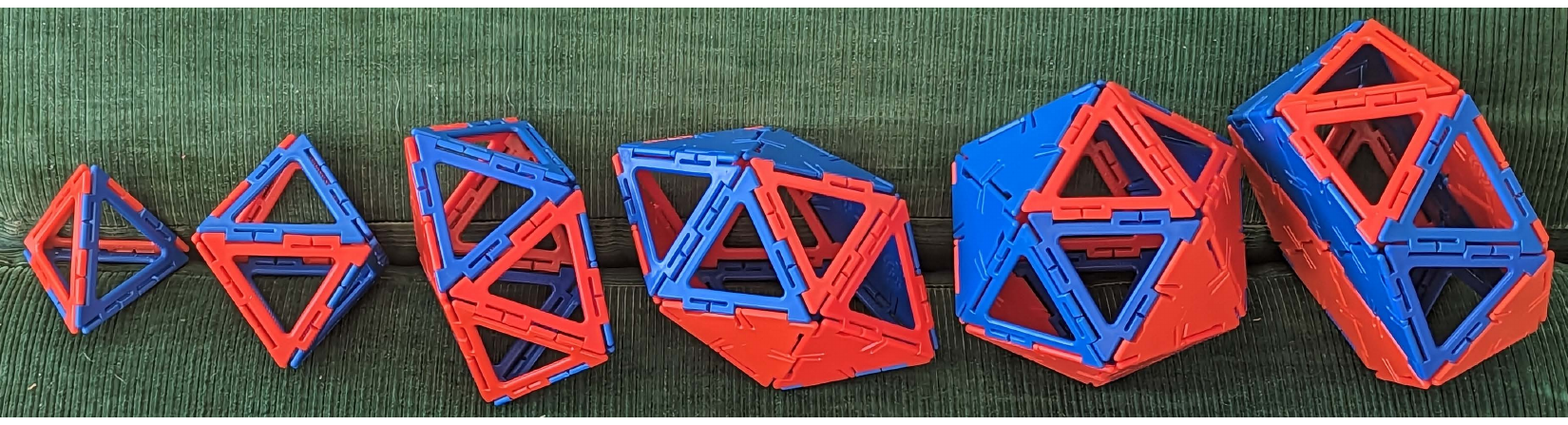}{icosplus}{Truncating $\Phi_{4,2}$
and capping in the
most natural way yields a family of phyllohedra $U_v$, $v$ even,
answering the question,
`What would the icosahedron look like if it had 14 vertices?'}
In fact there is a volume potential $\rho$ on the infinite cylinder which
when restricted may come in lower than the actual value of $\Qvol$,
but still yields $\tetvol$ when rounded up, giving
\[
\tetvol(U_v) = \ceiling{\rho(U_v)} = \ceiling{31/32(2v-12)+5/2}
.
\]
(See Figure \ref{defect}.)
\begin{figure}
\[
\begin{array}{r|r|r|r}
v&\tetvol(U_v)&\tetvol-\Qvol&\tetvol-\rho\\
\hline
12&15&0&7/8\\
14&18&0&0\\
16&22&0&1/8\\
18&26&0&1/4\\
20&30&1/5&3/8\\
22&34&1/3&1/2\\
24&38&1/2&5/8\\
26&42&2/3&3/4\\
28&46&2/3&7/8\\
30&49&0&0\\
32&53&0&1/8
\end{array}
\]
\caption{A single volume potential $\rho$ defined on $\Phi_{4,2}$
produces sharp lower bounds $\ceiling{\rho(U_v)}$ for $\tetvol(U_v)$.
The table shows the shortfall of the lower
bounds $\Qvol(U_v)$ and $\rho(U_v)$. As these are less than $1$,
rounding up gives the exact value.}
\label{defect}
\end{figure}

For $\Phi_{3,3}$ the results are similar, only now
with $\vsa= \frac{23}{24}$.

For girth 7 all $\vsa$ are $1$.
From $\Phi_{5,2}$ or $\Phi_{4,3}$ we get two new families with
$\tetvol=2v-10$.
For $\Phi_{6,1}$, when you cap in the natural way there are vertices of degree
$7$, and $\tetvol=2v-11$;
other methods of capping give triations with $\tetvol=2v-10$.

\begin{comment}
These truncated phyllocylinders that we've been discussing
can be realized in 3-space
as surfaces built out of equilateral triangles
(see Figures \ref{snub} and \ref{nelson}).
At least, the triangles are \emph{nearly} equilateral.
This brings to mind the story of Ed Nelson and the pentagon.
When Ed was in kindergarten, someone showed him a drawing of a
regular pentagon.
Ed's response was, `Interesting! I wonder if it exists.'
\fig{nelson}{nelson}{Do I exist?}
\end{comment}

\section{Flip distance} \label{flipdist}
The work of STT
\cite{stt:jams}
on tetrations
was motivated by the problem of finding the maximum possible
flip distance between two triangulations of a $v$-gon.
If $\alpha,\beta$ are triangulations of a $v$-gon with no common edge
we get a triation
$\alpha-\beta$
of the $2$-sphere by gluing along their common $v$-gon.
A flip path from $\alpha$ to $\beta$ begets a tetration of $\alpha-\beta$
so
\[
\flipdist(\alpha,\beta) \geq \tetvol(\alpha-\beta)
.
\]

Pournin
\cite{pournin:diameter}.
gave examples with $\flipdist=2v-10$, $v \geq 13$,
thus proving the analog of the tetvol conjecture for $\flipdist$.
We can identify his examples as arising
from truncations of $\Phi_{5,2}$,
outfitted with a particular natural Hamiltonian cycle.
By varying the Hamiltonian cycles
we get many other examples from these same triations.

Pournin's $\flipdist$ result doesn't imply the tetvol conjecture for
triations, because there may be a gap between $\tetvol$ and $\flipdist$:
The Hamiltonian cycle can prevent us from converting an optimal tetration
into a flip path.
Wang
\cite{wang:gap}
gives examples
where the ratio of $\flipdist$ to $\tetvol$
is arbitrarily close to $3/2$.
In Wang's examples $\flipdist$ and $\tetvol$ are on the order of $3/2v$ and
$v$, so far from the kind of triations we're dealing with here.
But we can take 
Wang's smallest example,
which has
$v=10$, $\flipdist=2v-10=10$, $\tetvol=2v-11=9$
(Figure \ref{v10}),
and inflate it to an example with
$v=16$,
$\tetvol(\alpha - \beta) = 2v-11 = 21$,
 and $\flipdist(\alpha,\beta)=2v-10 = 22$.
(Figures \ref{v16},\ref{snub}.)
This illustrates why knowing examples with $\flipdist=2v-10$
doesn't settle the tetvol conjecture.
\figsize{.8\linewidth}{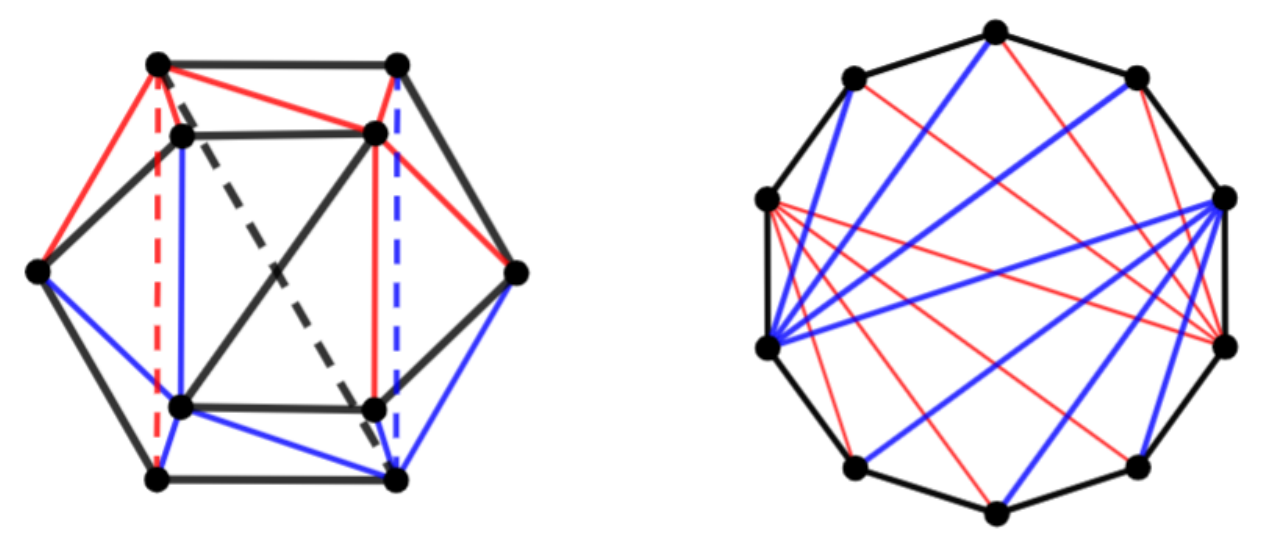}{v10}{The 10-vertex triation on the left has $\tetvol=2v-11=9$.
The black Hamiltonian cycle divides it into disk
triangulations $\alpha$ (blue) and $\beta$ (red),
with $\flipdist(\alpha,\beta)=2v-10 = 10$.}
\figsize{.8\linewidth}{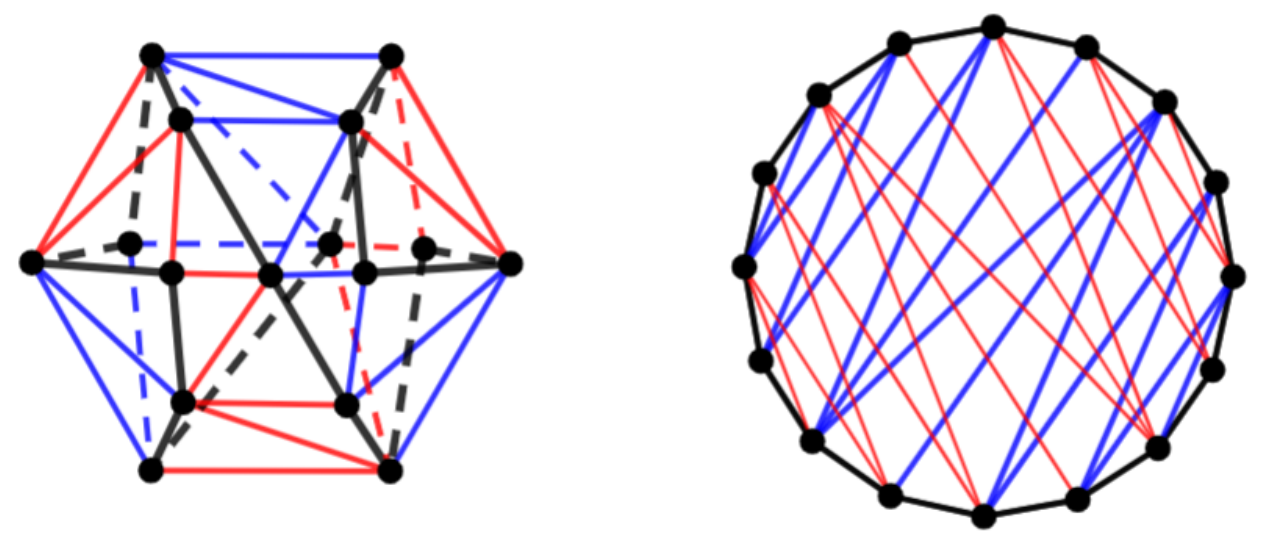}{v16}{A larger example.
This time $v=16$,
$\tetvol(\alpha - \beta) = 2v-11 = 21$,
and $\flipdist(\alpha,\beta)=2v-10 = 22$.}
\figsize{\linewidth}{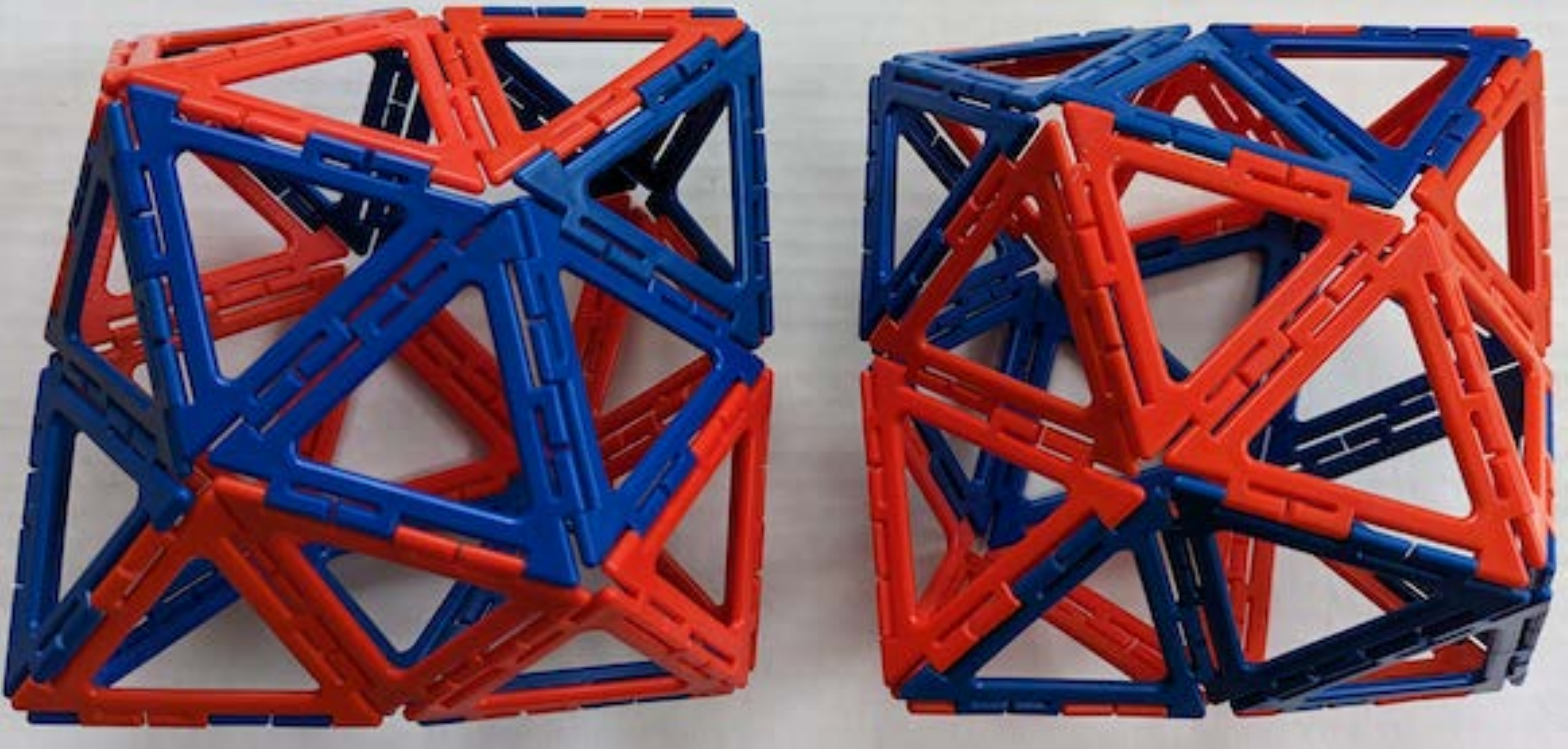}{snub}{The $v=16$ example brought to life.}

\section{The 3-ball}
We've sidestepped the question of whether a minimal tetration
of a triation $\sigma$ of the 2-sphere
necessarily yields a triangulation of the $3$-ball. It does.
That's because there are fewer tets than faces, so some tet must
meet $\sigma$ in at least two faces, necessarily adjacent.
Remove this tet, and you get a minimal tetration of a sphere, or a pair
of spheres joined along an edge
(or nothing at all, if you were down to a single tet).
Proceed by induction.
The key to the proof is the general
fact that any optimal filling of an integral $n$-cycle
splits under
what we call almost disjoint union,
where summands are supported on sets that
overlap in at most $n+1$ vertices.
See \cite{dew:taut} for details.

Contrast this behavior for triations of spheres
with what happens for triations of a torus,
or a surface of higher genus.
In that case there is no guarantee that a minimal tetration will be
a manifold, or even a pseudo-manifold.
And if it is a manifold,
we have no a priori control over how it fills in the surface.
All very mysterious.

\section{What's true in general}

If a triation has all vertices of degree 5 or 6, the chances are that
$\tetvol = 2v-10$.
The only exceptions we know are phyllohedra derived from $\Phi_{4,2}$ and
$\Phi_{3,3}$.
(This includes the icosahedron, which can be derived from either.)
We'll stop short of formulating a precise conjecture.
The point is that producing $2v-10$ triations is not the issue,
it's proving that they have this property.
As stated in section \ref{mt} above,
we expect that the right approach is that of MT.

\bibliography{bill}
\bibliographystyle{plain}
\end{document}